\documentclass[12pt]{article}

\usepackage{amsfonts,amsmath,amsthm,a4}
\usepackage{graphicx}
\usepackage{caption}
\usepackage{setspace}
\usepackage{natbib}

\usepackage{floatflt}

\newcommand{\Prob}{\mathbb{P}}
\newcommand{\N}{\mathbb{N}}	
\newcommand{\Z}{\mathbb{Z}}	
\newcommand{\R}{\mathbb{R}}	

\newcommand{\calX}{\mathcal{X}}
\newcommand{\calY}{\mathcal{Y}}
\newcommand{\Yhat}{\widehat{\mathcal{Y}}}
\newcommand{\Op}{\mathcal{F}}
\newcommand{\Ophat}{\widehat{\Op}}
\newcommand{\Fhat}{\widehat{\Op}}

\newcommand{\unknown}{\varphi}
\newcommand{\xdag}{\unknown^{\dagger}}

\newcommand{\its}{k}
\newcommand{\capit}{K}
\newcommand{\alk}{\alpha_{\its}}

\newcommand{\T}{\mathcal{T}}
\newcommand{\linOp}{\T}

\newcommand{\Tkm}{\widehat{\mathcal{T}}_{\its-1}}
\newcommand{\Tdag}{\widehat{\mathcal{T}}}
\newcommand{\hatxk}{\widehat{\unknown}_{\its}}
\newcommand{\hatxkm}{\widehat{\unknown}_{\its-1}}
\newcommand{\calR}{\mathcal{R}}
\newcommand{\ek}{e_{\its}}
\newcommand{\ekm}{e_{\its-1}}
\newcommand{\Delk}{\Delta_{\its}}
\newcommand{\Ex}{\mathbb{E}}

\newcommand{\lsp}{\left\langle}
\newcommand{\rsp}{\right\rangle}
\newcommand{\argmin}{\mathop{\mathrm{argmin}}}
\newcommand{\paren}[1]{\left( #1 \right)}

\def\Var{\mathop{\rm {\mathbb V}ar}\nolimits}%

\newcommand\independent{\protect\mathpalette{\protect\independenT}{\perp}}
\def\independenT#1#2{\mathrel{\rlap{$#1#2$}\mkern2mu{#1#2}}}

\newtheorem{theorem}{Theorem}
\newtheorem{corollary}[theorem]{Corollary}
\newtheorem{lemma}[theorem]{Lemma}

\begin{document}
\title
  {Iterative Estimation of Solutions to Noisy Nonlinear Operator Equations 
  in Nonparametric Instrumental Regression}
\author{Fabian DUNKER\footnote{Institute of Numerical and Applied Mathematics, University of G\"ottingen, Lotzestr. 16--18, 37083 G\"ottingen, Germany}
\footnote{Corresponding author \textit{Email}: dunker@math.uni-goettingen.de \textit{Tel.}: +49551394507} ,
Jean-Pierre FLORENS\footnote{Universite de Toulouse (GREMAQ and IDEI)} ,\\
Thorsten HOHAGE\footnote{Institute of Numerical and Applied Mathematics, University of G\"ottingen, Lotzestr. 16--18, 37083 G\"ottingen, Germany} , 
Jan JOHANNES\footnote{Institut de statistique, UCL, Voie du Roman Pays, 20, 1348 Louvain-la-Neuve Belgium} ,\\
Enno MAMMEN\footnote{Department of Economics, University Mannheim, L7,3-5, 68131 Mannheim, Germany
}}

\date{February 10, 2013}

\maketitle

\parindent0pt\textbf{Abstract:} This paper discusses the solution of 
nonlinear integral equations with noisy integral kernels as they appear
in nonparametric instrumental regression. We propose a regularized
Newton-type iteration and establish convergence and convergence rate results. 
A particular emphasis is on instrumental regression models where the usual
conditional mean assumption is replaced by a stronger independence assumption. 
We demonstrate for the case of a binary instrument that our approach allows
the correct estimation of regression functions which are not identifiable 
with the standard model. This is illustrated in computed examples with 
simulated data. 

 \vskip .2in
\noindent {\sl MSC: AMS 2000 subject classification.}
primary 62G08, secondary 62G20
\newpage

\onehalfspacing

\section {Introduction}

In this paper we will propose and analyze an iterative method for estimating
the solution of nonlinear integral equations which appear in nonparametric instrumental
regression problems. Examples will be discussed below, see  eq.~\eqref{eq:instr_model} and Section \ref{sec:example}. 
Such integral equations can be written as nonlinear operator equations
\begin{equation}\label{eq:nonlinop}
\Op(\varphi)=0
\end{equation}
where the operator $\Op$ is unknown, but where an estimator 
$\Fhat$ of $\Op$ is available. We will assume that $\Op : \mathfrak{B} \subset \calX \rightarrow \calY$
maps from a convex set $\mathfrak{B}$ in a Banach space $\calX$ to a Hilbert 
space $\calY$. Typically such operator equations are ill-posed in the sense 
that $\Op^{-1}$ is not continuous. In particular this is the case for
integral operators with smooth kernels on a compact set. 
In such cases the straightforward estimator 
$\Fhat^{-1}(0)$ will not be consistent since it has infinite variance. 
Regularization techniques must be applied to
solve \eqref{eq:nonlinop} or its empirical version $\Fhat(\hat \varphi)=0$.
We will use a generalized version of
the iteratively regularized Gau{\ss}-Newton
method. In numerical analysis this is one of the most popular computational methods for solving nonlinear
ill-posed operator equations.
It avoids some problems of  nonlinear Tikhonov regularization given by
\begin{equation}\label{eq:nltikh}
\widehat{\varphi}:=\argmin_{\unknown}
\left[\|\Fhat(\unknown)\|_\calY^2 + \alpha \|\unknown-\unknown_0\|_\calX^2\right],
\end{equation}
where  $\unknown_0$ is some initial guess of $\varphi$.
In practice the iteratively regularized Gau{\ss}-Newton
method does not suffer from the problem that minima of
the functional in \eqref{eq:nltikh} need not to be unique and it avoids computational difficulties
due to the presence of local minima. We will compare both methods in more details later.
Moreover, instead of a quadratic penalty,
we allow for a more general penalty term $\mathcal{R}:\mathfrak{B}\to (-\infty,\infty]$ with domain of
definition $\mathfrak{B}$. We only assume that  $\mathcal{R}$ is a convex, lower semi-continuous
functional  that is  not identically equal to $\infty$.
With this choice an iteratively regularized Gau{\ss}-Newton
method is given by the iterations
\begin{equation}\label{eq:defi_method}
\hatxk := \argmin_{\unknown\in \mathcal{B}}
\left[\|\Fhat'[\hatxkm](\unknown-\hatxkm) + \Fhat(\hatxkm)\|_\calY^2 + \alk \calR(\unknown) \right].
\end{equation}
In each Newton step a convex optimization problem has to be solved with a sequence
of regularization parameters $\alk$. We assume that $\alk$ tends to $0$ in a way that will be specified in Section
\ref{sec:convergence}. In the special case that $\calX$ is a Hilbert space, the most common choice for the penalty term is 
$\calR(\unknown) = \|\unknown-\unknown_0\|_{\mathcal{X}}^2$. Here $\|\cdot\|_{\mathcal{X}}$
is the norm of the Hilbert space $\mathcal{X}$ and $\unknown_0$ is the initial guess
at which the iteration is started. This is the iteratively regularized Gau{\ss}-Newton
method as suggested by \cite{bakush:92} and further analyzed by \cite{BNS:97} and \cite{hohage:97} 
for low order H\"older or logarithmic source conditions, respectively.
We also refer to the monographs by \cite{BK:04b} and \cite{KNS:08} and to further references therein.

The use of more general convex regularization terms in the general case where $\calX$ is a Banach space allows for a flexible
incorporation of further a-priori information. Common choices are entropy regularization,  $l^1$ penalties and bounded 
variation ($\mathrm{BV}$) penalties. \cite{LouPel:08} studied entropy regularization for instrumental 
variable models but they gave no theoretical results for the rates of convergence of their estimators. If a basis or a 
frame of $\calX$ is given, an $l^1$ penalty of the coefficients with 
respect to this basis or frame enhances sparsity properties of the estimator 
with respect to this basis or frame.  A $\mathrm{BV}$ penalty is particularly
appropriate for piecewise constant estimates.

Our main result gives rates of convergence for the estimator where the distance between the estimator and the solution of
\eqref{eq:nonlinop} is measured by the Bregman distance, see Theorem \ref{theo:main}.
For entropy regularization this directly implies  convergence estimates measured by the $L^1$-norm. Our scheme allows for
the incorporation of structural a-priori information of the form
$\unknown\in \mathcal{C}$ where   $\mathcal{C}$ is a closed convex set
(e.g.~a-priori information on non-negativity, monotonicity or convexity/concavity). This can be done by
setting $\mathcal{R}(\unknown):=\infty$ if $\unknown \notin \mathcal{C}$. 

For convex regularization terms, the analysis differs from the mathematical approaches used for studying quadratic
regularization. One has to employ variational methods
rather than spectral methods. Recently, a number of papers have appeared on this subject, we only mention \cite{eggermont:93}, \cite{BO:04}, \cite{resmerita:05}, \cite{HKPS:08}, \cite{scherzer_etal:09}.
A first variational convergence rate analysis
of Newton-type methods in a deterministic setting without errors in the operator
and $\mathcal{R}$ given by Banach norms has recently been done by \cite{KH:10}.
Our analysis is closest to that of the last reference. However, all the
references above only treat perturbations of the right hand side of 
the operator equation, and hence these results are not applicable to 
nonparametric instrumental regression. Our treatment of nonlinear ill-posed 
operator equations with errors in the operator may be of independent
interest and relevant for other applications. 

For the special case that $\calX$ is a Hilbert space convergence rates of the nonlinear Tikhonov regularization were discussed 
in \cite{EKN:89} in a deterministic setting. Rates for a model with random errors
were obtained in \cite{BisHohMun:04}. In \cite{HorLee:07} nonparametric instrumental variables estimation is 
considered in a quantile regression model. This is one example of a statistical model where the unknown nonparametric function 
is given as the solution of a nonlinear integral equation. We  will describe this model in the next section. 

In \cite{HorLee:07} it is assumed that the singular values of the Fr\'echet derivative $\Op'[\unknown]$ decay polynomially
and results are given on the rates of convergence under these assumptions. Horowitz and Lee pointed out that a convergence
analysis for exponentially decreasing singular values is an important open problem. We will show that singular values of 
integral operators with infinitely smooth kernels do in fact decrease super-algebraically
and present a convergence analysis without an assumption on the rate of decay of the singular values.

Besides the analysis of the iteratively regularized Gau{\ss}-Newton method for noisy operators the second main innovation of this
paper is a nonparametric instrumental regression models where the instrument $W$ is independent from the error $U$:
\begin{subequations}\label{eq:instr_model}
\begin{eqnarray}
\label{eq:indep1}
&&Y = \varphi(Z)+U,\\
\label{eq:indep_ass}
&&U \independent W,\\
\label{eq:meanU}
&&\mathbb{E}U =0.
\end{eqnarray}
\end{subequations}
Here, $Y$ is a scalar response variable, $Z$ is an observed random vector of endogeneous explanatory variables.
It is shown in Section \ref{sec:example} that this model leads to a nonlinear
integral equation of the form \eqref{eq:nonlinop}
with a kernel, that has to be estimated from data.

This model slightly differs from nonparametric instrumental regression with mean independent instruments given by
\begin{subequations}\label{eq:standard_instr}
\begin{eqnarray}
&&Y = \varphi(Z) + U,\\
\label{eq:conditional_mean}
&&\mathbb{E}[U|W] = 0.
\end{eqnarray}
\end{subequations}
The latter model has been studied intensively in econometrics by a number of authors, see e.g. \cite{florens:03},
\cite{NewPow:03}, \cite{HalHor:05}, \cite{BluCheKri:07}, \cite{CheRei:10} and \cite{BreJoh:09}.
In this model the regression function $\varphi$ is defined as the solution of a linear first kind integral equation
\begin{equation}\label{eq:linop}
\mathcal{T}\varphi=g
\end{equation}
where both the kernel of the integral operator
$(\mathcal{T}\varphi)(w):=\mathbb{E}[\varphi(Z)|W=w]$ and the right hand side $g(w):= \mathbb{E}[Y|W=w]$ have to be estimated from the data.

Actually, in specific econometric applications, the conditional mean assumption \eqref{eq:conditional_mean} is typically established
by arguing that the stronger
independence assumption \eqref{eq:indep_ass} holds. Therefore, it is a natural question if one can improve the accuracy of estimation of $\varphi$  by using
the stronger condition \eqref{eq:meanU}, \eqref{eq:indep_ass} directly. We will give a first partial positive answer to this question: a necessary condition 
for identifiability in the model \eqref{eq:standard_instr} is that the instrumental variable $W$ must have at least as many continuously distributed 
components as the explanatory variable $Z$. This is not necessary in model  \eqref{eq:instr_model}. As an example we will demonstrate that $\varphi$ can be 
identifiable even if $W$ is binary and $Z$ is one-dimensional and continuously distributed. Hence, the model \eqref{eq:instr_model}
contains strictly more information on $\varphi$ than the model \eqref{eq:standard_instr}.
A more detailed comparison of the two models is very complex because  the integral equations
obtained from these two models are related only very implicitly.

The plan of this paper is as follows: in the following section
we give more details on our motivating examples from instrumental variable regression.
Section \ref{sec:sc} recalls the definition of source conditions and discusses their relation to smoothness conditions. 
In particular, we show that
for integral equations of the first kind with smooth kernels, H\"older type
source conditions are too restrictive, and discuss variational forms of
source conditions. In Section \ref{sec:convergence} we present our main convergence result for
the iteratively regularized Gau{\ss}-Newton method with noisy operators. Afterwards, we discuss in Section \ref{sec:revisited} 
how this result applies to the regression problem \eqref{eq:instr_model}. Section \ref{sec:simulations} reports on numerical
simulations for an instrumental variable regression model with binary instruments.

\section{Examples}\label{sec:example}
\subsection{Instrumental quantile regression}\label{sec:example:quantile}

In \cite{HorLee:07}  the following quantile regression model has been studied: 
\begin{subequations}\label{eq:quantile_regression}
\begin{eqnarray}
&&Y= \varphi(Z) +U \\
&& \Prob(U\leq 0|W=w) = q \qquad \mbox{for all }w
\end{eqnarray}
\end{subequations}
Here,  $Y$ is a response variable, $Z$ is an endogeneous explanatory variable, $q\in (0,1)$ is a fixed constant,
$U$ is an unobserved error variable and $W$ an observable instrument. The quantile is defined conditional on $W$.

We assume from now on that each of the random variables $Y$, $Z$ and $W$ is a vector of continuous or discrete random variables. Further, we assume that a joint density $f_{YZW}$ exists with respect to the
Lebesgue measure, the counting measure or a product of both measures respectively. 
Let $G_{YZW}(y,z,w):=\int_{-\infty}^y f_{YZW}(\tilde{y},z,w)\,d\tilde{y}$,
and let $f_W(w):=\int\int f(y,z,w)\,dy\,dz$ denote the marginal density of $W$.
Then $\varphi$ solves a nonlinear operator equation \eqref{eq:nonlinop} with the operator
\[
(\Op(\varphi))(w) := \int G_{YZW}(\varphi(z),z,w)\,dz- q f_W(w).
\]
It is pointed out in \cite{HorLee:07}, that the model \eqref{eq:quantile_regression} subsumes
nonseparable quantile regression models of the form
\begin{equation}\label{eq:nonsep_quantile_regr}
Y=H(Z,V)
\end{equation}
as studied in Chernozhukov \& Imbens \& Newey \cite{CheImbNew:07}, see also Chernozhukov \& Hansen \cite{CheHan:05}. Here $V$ is an unobserved,
continuously distributed random variable independent of an instrument $W$, and the function
$H$ is strictly increasing in its second argument. Assuming w.l.o.g.\ that
$V\sim U[0,1]$, \eqref{eq:nonsep_quantile_regr} reduces to \eqref{eq:quantile_regression}
with $U:=Y-H(Z,q)$ and $\varphi(z):=H(z,q)$.

\subsection{Nonparametric regression with independent instruments}\label{sec:example:indep}

\subsubsection{Operator equations}
The model \eqref{eq:indep1}, \eqref{eq:indep_ass} leads to the nonlinear integral equation
\begin{subequations}\label{eqs:indep_ieq}
\begin{equation}\label{eq:inteq}
\int f_{YZW}(u+\varphi(z),z,w)\,dz -
\int f_{YZ}(u+\varphi(z),z) f_W(w)\,dz=0,\quad
\mbox{for all }u, w,
\end{equation}
where we assume as above that the joint density $f_{YZW}$ of $(Y,Z,W)$ exists. The marginal densities
of $(Y,Z)$ and $W$ are denoted by $f_{YZ}$ and $f_{W}$ respectively. Note that if $\varphi$ is
a solution to \eqref{eq:inteq}, then any function $\varphi+a$ with $a\in\R$
is another solution to \eqref{eq:inteq}. The additive constant can be fixed
by taking into account eq.~\eqref{eq:meanU}, which may be rewritten as
\begin{equation}\label{eq:scale}
\int \varphi(z)f_Z(z)\,dz-\int y f_Y(y)\,dy  = 0
\end{equation}
\end{subequations}
with the marginal densities $f_Y$ and $f_Z$ of $Y$ and $Z$.

The system of equations \eqref{eq:inteq}, \eqref{eq:scale} 
can be written as a nonlinear ill-posed operator equation
\eqref{eq:nonlinop} with the operator
\begin{equation}\label{eq:nonlininstreg:opeq}
(\Op(\varphi))(u,w) \! := \! \left(\begin{array}{c}
\int f_{YZW}(u+\varphi(z),z,w) - f_{YZ}(u+\varphi(z),z) f_W(w)\,dz \\
\int \varphi(z)f_Z(z)\,dz-\int y f_Y(y)\,dy
\end{array}\right) \! .
\end{equation}
If we assume the existence of the joint density of $(Y,Z)$ unconditional and conditional
given $W$, say $f_{Y,Z}$ and $f_{Y,Z|W}$ respectively we can use the equivalent operator
\begin{equation}\label{eq:nonlininstreg:opeq3}
(\bar\Op(\varphi))(u,w) := \left(\begin{array}{c}
\int f_{YZ|W}(u+\varphi(z),z|w) - f_{YZ}(u+\varphi(z),z)\,dz \\
\int \varphi(z)f_Z(z)\,dz-\int y f_Y(y)\,dy
\end{array}\right) .
\end{equation}
Alternatively, it may be advantageous to integrate \eqref{eq:nonlininstreg:opeq3} once with respect to $u$.
Introducing $G_{YZ|W}(\tilde{y},z|w):=\int_{-\infty}^{\tilde{y}}f_{YZ|W}(y,z|w)\,dy$
and $G_{YZ}(\tilde{y},z):=\int_{-\infty}^{\tilde{y}} f_{YZ}(y,z)\,dy$ yields
an other operator formulation of the model \eqref{eq:instr_model} with the operator
\begin{equation}\label{eq:nonlininstreg:opeq2}
(\tilde{\Op}(\varphi))(\tilde{u},w) := \left(\begin{array}{c}
\int G_{YZ|W}(\tilde{u}+\varphi(z),z|w) - G_{YZ}(\tilde{u}+\varphi(z),z) \,dz \\
\int \varphi(z)f_Z(z)\,dz-\int y f_Y(y)\,dy
\end{array}\right).
\end{equation}
Let us set $(\mathcal{G}(\varphi))(u,w) := \int G_{YZ|W}(\tilde{u}+\varphi(z),z|w) - G_{YZ}(\tilde{u}+\varphi(z),z) \,dz$. Then the last operators can be written as
\[
\tilde\Op(\varphi):=\left(\begin{array}{c}
\mathcal{G}(\varphi)\\
\Ex(Y-\varphi(Z))
\end{array}\right).
\]

\subsubsection{Identification}
In  the following we discuss sufficient conditions for the injectivity of the Gateaux derivative $\tilde\Op'[\varphi]$ of $\tilde\Op$
at the
solution $\varphi$. Local identifiability of the nonlinear problem $\mathcal{F}(\varphi)=0$ in an open neighborhood
of $\varphi$  is not necessarily implied by injectivity of  $\mathcal{F}'[\varphi]$. Additional assumptions that guarantee local identifiability are Frechet differentiability and 
 tangential cone conditions, compare (\ref{eq:tcn}). For a discussion we refer to  \cite{KNS:08}, \cite{CCLN:11} or \cite{FloSba:10}.
The Gateaux derivative of $\tilde\Op$ at $\varphi$ is given by 
\begin{align*}
\tilde\Op'[\varphi]\phi = \begin{pmatrix} \mathcal{G}'[\varphi]\phi\\\Ex(\phi(Z))\end{pmatrix}
\end{align*}
where the Gateaux derivative $\mathcal{G}'[\varphi]$ of $\mathcal{G}$ at $\varphi$ satisfies
\begin{align}\label{eq:nonlininstreg:derivative}
(\mathcal{G}'[\varphi]\phi)(u,w) = \int \phi(z) \big( f_{Y,Z|W}(u+\varphi(z),z|w)- f_{Y,Z}(u+\varphi(z),z)\big)\, dz.
\end{align}
Injectivity of $\tilde\Op'[\varphi]$ is equivalent to injectivity of $\mathcal{G}'[\varphi]$ on the linear subspace of functions $\phi$ with 
$\Ex[\phi(Z)]=0$. We denote by $f_{U}$ the marginal density of $U = Y-\varphi(Z)$. Then  by employing the independence of $U$ and $W$ a 
change of variables allows us to write
\begin{equation}
\mathcal{G}'[\varphi]\phi=\bigg(\Ex[\phi(Z)|U,W]-\Ex[\phi(Z)|U]\bigg) f_U.
\end{equation}
Alternatively, we may consider the linear operator
\begin{equation}\label{eq:3}
\linOp\phi:= \Ex[\phi(Z)|U,W]-\Ex[\phi(Z)|U]
\end{equation}
mapping from a function space of mean zero functions in $Z$ into a function space in $U \times W$. Roughly speaking, injectivity of the operator $\linOp$ and hence local identification is possible  if the dependence between the endogenous regressor $Z$ and the error term $U$ varies sufficiently with
respect to the instrument $W$. The next example illustrates this fact.

\paragraph{Example 1.} Let $U$, $V$ and $W$ be  real valued independent random variables and  let $\rho$ be a function
defined on $\R$ and taking values in $[-1,1]$. Define the endogenous regressor
\[
Z := U\,\rho(W) +V\,\sqrt{|1-\rho^2(W)|}.
\]
If $U$ and $V$ are standard normally distributed, which we assume in this example, then it is easily seen that the conditional distribution
of $(U,Z)$ given $W$ is Gaussian:
\begin{equation}
\begin{pmatrix}U\\Z\end{pmatrix}\bigg|W\sim \mathcal{N}\left(\begin{pmatrix}0\\0\end{pmatrix},\begin{pmatrix}1&\rho(W)\\\rho(W)&1\end{pmatrix}\right).\label{eq:4}
\end{equation}
Note that in this situation  $U$ and $Z$ are marginally standard normally distributed, both  unconditional and conditional on $W$.
In other words, $U$ and $W$ as well as $Z$ and $W$ are independent. But obviously, the random vector $(U,Z)$ and the instrument  $W$
are dependent. Interestingly, in the commonly studied case of mean independence, that is $\Ex[U|W]=0$, identification is guaranteed  if
and only if the conditional distribution of $Z$ given $W$ is complete (cf. \cite{CFR06handbook}) which rules out  the independence of $Z$
and $W$ and hence this example. However, in this example the linear operator $\linOp$ defined in \eqref{eq:3} can  be injective and thus
local  identification might be still possible. In order to provide sufficient conditions to ensure injectivity of $\linOp$, let us recall
the eigenvalue decomposition of the conditional expectation operator for normally distributed random variables.

The following development can be found in \cite{CFR06handbook} while it has been shown  thoroughly in \cite{Letac1995}. Consider random
variables $U^*$ and $Z^*$ satisfying 
\[
\begin{pmatrix}U^*\\Z^*\end{pmatrix} \sim \mathcal{N} \left( \begin{pmatrix}0\\0\end{pmatrix}, \begin{pmatrix}1&\rho\\\rho&1\end{pmatrix} 
\right)
\] 
for some $\rho\in[-1,1]$. Obviously, $U^*$ and $Z^*$ are marginally identically distributed with  standard normal density $f_{0,1}$ which
in turn implies $L^2_{U^*}=L^2_{Z^*}=:L^2_{f_{0,1}}$. Note, that by an elementary symmetry argument the  conditional expectation operator
$S\phi:=\Ex[\phi(Z^*)|U^*]$ of $Z^*$ given $U^*$ mapping $L^2_{f_{0,1}}$ to itself is self-adjoint and hence  $S$ permits an eigenvalue
decomposition. Moreover, for $j=0,1,2,\dotsc$ let $f^{(j)}_{0,1}$ be the $j$th derivative of $f_{0,1}$ and  let 
$H_j:=(-1)^jf^{(j)}_{0,1}/f_{0,1}$ denote the $j$th Hermite polynomial.  The Hermite polynomials form a complete orthogonal system in $L^2_{f_{0,1}}$, see e.g.\ Problem IV-29 on page 117 in \cite{Letac1995}. 
Furthermore,  $\Ex[H_j(Z^*)|U^*]=H_j(U^*)\rho^j$ holds true for all $j \in \N_0$, see e.g.\  Problem IV-30 on page 120 in \cite{Letac1995}. From these assertions we readily conclude that the
eigenfunctions $\{\psi_j\}_{j=0}^\infty$ of $S$  are  up to multiples given by the Hermite polynomials and that 
$(\rho^j)_{j \in \N_0}$ is the corresponding sequence of eigenvalues.

Keeping in mind  that the distribution of $(U,Z)$ conditional on $W$ given in \eqref{eq:4} is Gaussian let us reconsider the operator
$\linOp$ defined in \eqref{eq:3}. By employing that $U$ and $W$ are independent it is straightforward to conclude that 
\[
\Ex[|(\linOp\phi)(U,W)|^2]=\sum_{j=1}^\infty\Var(\rho^j(W))\Ex[|\phi(Z)\psi_j(Z)|^2]
\]
for all $\phi\in L^2_Z$ with $\Ex[\phi(Z)]=0$, where the basis $\{\psi_j\}_{j=1}^\infty$  are   multiples of the Hermite polynomials. Consequently,  the operator $\linOp$ is injective if
and only if  $\Var(\rho^j(W))\ne0$ for all $j\in\N$, (keep in mind Parseval's identity, i.e. 
$\Ex[f(Z) ^2]=\sum_{j=1}^\infty\Ex[f(Z)\psi_j(Z)]^2$ for all $f\in L^2_Z$). This in turn holds if and only if the random variable 
$|\rho(W)|$ is not constant. Surprisingly, even in case of a binary instrument $W$ taking only two values, say
$P(W=0)=w_0$ and $P(W=1)=1-w_0$ with $0<w_0<1$, the condition $|\rho(0)|\ne|\rho(1)|$ is sufficient to ensure  the injectivity of
the operator $\linOp$. 

\paragraph{Example 2.} We now give another example for injectivity of $\linOp$. We consider again a binary instrument $W$ and we make the additional assumption that the conditional copula function of $U$ and $Z$, given $W=w$ does not depend on $w$. This assumption has been made by \cite{ImbensNewey2009} in case of a continuous instrument. Under this assumption it holds that $(U,V)$ is independent of $W$ where $V=F_{Z|W}(Z|W)$ for the conditional distribution function $F_{Z|W}(z|w)$ of $Z$ given $W=w$. Note that in the case of a binary instrument injectivity of $\linOp$ is equivalent to the injectivity of the map
\[
\phi \mapsto  \Ex[\phi(Z)|U, W=1] -   \Ex[\phi(Z)|U, W=0]
\]
on the space of all functions $\phi$ with $ \Ex[\phi(Z)]=0$. We use that 
\begin{eqnarray*}
0& =&\Ex[\phi(Z)|U, W=1] -   \Ex[\phi(Z)|U, W=0] \\ &=& \Ex[\phi(F^{-1}_{V|W}(V|1))|U, W=1] -   \Ex[\phi(F^{-1}_{V|W}(V|0))|U, W=0]  \\ &=& \Ex[\phi(F^{-1}_{V|W}(V|1))|U] -   \Ex[\phi(F^{-1}_{V|W}(V|0))|U] \\ &=&  \Ex[\phi(F^{-1}_{V|W}(V|1))-\phi(F^{-1}_{V|W}(V|0))|U], 
\end{eqnarray*}
because of independence of $(U,V)$ and $W$. If the family of conditional densities of $V$ given $U$ is complete this equation implies that $\phi(F^{-1}_{V|W}(v|1))=\phi(F^{-1}_{V|W}(v|0))$ almost surely. The latter equation can be used to get that under some additional assumptions on $F_{Z|W}$ the function $\phi$ is almost surely constant, see the arguments used in \cite{Torgo} and \cite{Hault}. Because of $ \Ex[\phi(Z)]=0$ we get that  $\phi(z) = 0$ a.s. Thus $\linOp$ is invertible. Note that our discussion differs from the results in \cite{ImbensNewey2009}, \cite{Torgo} and \cite{Hault}. We make the assumption on the conditional copula function only  for the underlying distribution and argue that - under additional conditions - local identifiability holds for a neighborhood of distributions for which this assumption may not apply whereas in the latter papers the conditional copula assumption is used as a model assumption for all distributions of the statistical model. This heuristic discussion can be 
generalized to more general instruments with discrete and/or continuous components.

We will continue the discussion of binary instruments in the next subsection. Section  \ref{sec:simulations} contains further numerical evidence of identifiability in a particular case.

\subsubsection{Binary instruments}

We consider the above mentioned special case that the instrument $W$ is binary and it only takes the values $0$ and $1$.
Furthermore, the explanatory variable $Z$ is a scalar. Then the marginal density $f_W$ (w.r.t. the counting measure) has the two values
\[
f_W(0)=w_0 \qquad \text{ and }\qquad  f_W(1)=w_1 = 1-w_0\,.
\]
Equation (\ref{eq:inteq}) is equivalent to the system of equations
\begin{equation*}
\begin{array}{l}
\int f_{YZW}(u + \varphi(z), z, 0)dz = w_0 \int f_{YZ}(u + \varphi(z), z)\,dz \\ [0.5ex]
\int f_{YZW}(u + \varphi(z), z, 1)dz = w_1 \int f_{YZ}(u + \varphi(z), z)\,dz
\end{array}
\qquad \text{for all } u\,.
\end{equation*}
It follows from the identity $f_{YZ}(y,z) = f_{YZW}(y,z,0) + f_{YZW}(y,z,1)$
that these two equations are linearly dependent and can be rewritten as
\begin{equation}\label{nullstelle}
\int w_1 f_{YZW}(u + \varphi(z), z, 0)
- w_0 f_{YZW}(u + \varphi(z), z, 1)\,dz = 0\quad \text{for all } u\,.
\end{equation}
So $\varphi$ is a root of the nonlinear ill-posed operator
\begin{equation}\label{eq:nonlinbininstreg:opeq}
(\Op(\varphi))(u) \! := \! \left(\begin{array}{c}
\int w_1 f_{YZW}(u + \varphi(z), z, 0)
- w_0 f_{YZW}(u + \varphi(z), z, 1)\,dz\\
\int \varphi(z)f_Z(z)\,dz-\int y f_Y(y)\,dy
\end{array}\right) \! .
\end{equation}
In analogy to \eqref{eq:nonlininstreg:opeq2}, the equation ${\Op}(\varphi)=0$
can equivalently be rewritten as $\tilde{\Op}(\varphi)=0$ with   
\begin{equation}\label{eq:nonlinbininstreg2:opeq}
(\tilde{\Op}(\varphi))(u) := \left( \! \begin{array}{c}
\int w_1 G_{YZW}(u + \varphi(z), z, 0)
- w_0 G_{YZW}(u + \varphi(z),z,1)\,dz\\
\int \varphi(z)f_Z(z)\,dz-\int y f_Y(y)\,dy
\end{array} \! \right) \!.
\end{equation}

We emphasize that $Z$ does not have to be discrete for identifiability, as it is the case when the conditional mean assumption \eqref{eq:conditional_mean} is used instead of the independence assumption \eqref{eq:meanU}. We will return to this point in Section \ref{sec:simulations}.

\section{Smoothness in terms of source conditions}\label{sec:sc}

In this section we collect some material on source conditions that will be needed 
in the next section to state our main result. We are primarily interested in source conditions 
in Banach spaces. However, we start with a motivation for $L^2$ spaces and present in a first 
step a definition of source conditions in the special case of Hilbert spaces. For the sake of 
simplicity we discuss the relevance of source conditions for nonparametric instrumental 
regression problems in this special case. Afterward, we introduce source conditions 
for the general case of Banach spaces.

\subsection{Source conditions in Hilbert spaces}

Let us recall the relationship between the smoothness of a kernel $k$ of a compact linear
integral operator $\linOp:L^2([0,1]^{d_1})\to L^2([0,1]^{d_2})$,
\[
(\linOp\varphi)(x) := \int_{[0,1]^{d_1}} k(x,y)\varphi(y)\, dy,
\qquad x\in [0,1]^{d_2}
\]
and the decay of its singular values $\sigma_j$.
If $\{(u_j,v_j,\sigma_j):j\in\mathbb{N}_0\}$ is a singular
system of $\linOp$, then according to the Courant-Fischer characterization
(see e.g.\ \cite{kress:99})
of the singular values the operator $\linOp_j$ with kernel
$k_j(x,y):=\sum_{l=0}^{j-1} \sigma_l v_l(x)u_l(y)$ satisfies
\begin{equation*}
\sigma_j=\|\linOp-\linOp_j\| = \inf\{ \|\linOp-\tilde{\linOp}\|: \mathrm{rank}\, \tilde{\linOp} \leq j\}.
\end{equation*}
In particular, if there exist functions
$\tilde{u}_l\in L^2([0,1]^{d_1})$,
$\tilde{v}_l\in L^2([0,1]^{d_2})$,
and numbers $\tilde{\sigma}_l$ for all $l\in \mathbb{N}_0$
such that $\int_{[0,1]^{d_1}}\int_{[0,1]^{d_2}}|k(x,y)-\sum_{l=0}^{j-1} \tilde{u}_l(x)\tilde{v}_l(y)|^2\,
dx\,dy\leq \tilde{\sigma}_{l}$, then $\sigma_j\leq \tilde{\sigma}_{j}$ since
$\|\linOp-\linOp_j\|\leq \|k-k_j\|_{L^2([0,1]^{d_1+d_2})}$. It follows from standard
results in approximation theory (see e.g.\ \cite{PS:91})
that for smooth bounded domains the singular values $\sigma_j$ decay
at least polynomially if $k$ belongs
to a Sobolev space, super-algebraically if $k\in C^{\infty}([0,1]^{d_1+d_2})$,
and at least exponentially if $k$ is analytic.

In regularization theory, smoothness of the solution $\xdag$
to an inverse problem is usually formulated in terms of source conditions, which describe
smoothness relative to the smoothing properties of the operator.
For a linear operator $\linOp:\mathcal{X}\to\mathcal{Y}$
between Hilbert spaces $\mathcal{X}$ and $\mathcal{Y}$,
such source conditions have the form
\begin{equation}\label{eq:classical_sc}
\xdag-\unknown_0 = \Lambda(\linOp^*\linOp)\psi\,.
\end{equation}
Here $\psi\in\mathcal{X}$, $\varphi_0$ is an initial guess (typically $\varphi_0=0$ in the linear case), $\linOp^*$ is the
adjoint operator of $\linOp$ with respect to the scalar product of the Hilbert space, and $\Lambda:[0,\infty)\to [0,\infty)$ is a continuous, strictly monotonically increasing function with $\Lambda(0)=0$.
$\Lambda(\linOp^*\linOp)$ is defined by using the spectral calculus. So
with the notations above $\Lambda(\linOp^*\linOp)\psi
=\sum_{l=0}^{\infty}\Lambda(\sigma_l^2)\langle \psi,u_l\rangle u_l$.
For a nonlinear operator between Hilbert spaces $\mathcal{F}:\mathcal{X}\to\mathcal{Y}$ the
Gateaux derivative $\linOp=\mathcal{F}'[\xdag]$ at $\xdag$
is used. 

If we choose a fixed $\Lambda$ the source condition is the more restrictive 
the faster the singular values decay. I.e. for integral operators it is the more restrictive
the smoother the kernel. For the most common choice $\Lambda(t)=t^{\mu}$ for some $\mu>0$
these condition are called \emph{H\"older-type source conditions}.
We refer to the monographs \cite{EHN:96,BK:04b,KNS:08} for further
information.

\subsection{Impact on nonparametric instrumental regression}\label{sec:source:impact}

Let us discuss source conditions in the context of nonparametric instrumental variable models. 
The kernel of the integral operator in \eqref{eq:nonlininstreg:derivative}
is composed of probability densities. For the derivatives of the alternative operators 
\eqref{eq:nonlininstreg:opeq} and \eqref{eq:nonlininstreg:opeq3} it is composed of partial derivatives of densities. Many typical probability density functions are analytic, i.e. the density of the normal. Hence, in applications it
will frequently occur that the kernel of the operator in the source condition is infinitely smooth or even analytic.

Let us have a closer look at these cases. The singular values of the operator in \eqref{eq:nonlininstreg:derivative} 
will decay super-algebraically or even exponentially. As a consequence, H\"older-type source conditions are extremely restrictive
smoothness conditions, since the eigenvalues $\lambda_j((\linOp^*\linOp)^{\nu}) = \sigma_j^{2\nu}$ will decay
super-algebraically or exponentially, too.  Hence, H\"older-type source conditions imply that
the Fourier coefficients with respect to $\{u_j:j\in\mathbb{N}_0\}$ of the difference between initial guess and regression function $\xdag-\unknown_0$ decay super-algebraically 
or exponentially. For standard Fourier coefficients this entails that $\xdag-\unknown_0$ has to be 
infinitely smooth or even analytic. Hence, the initial guess must be very good and already capture some features of the unknown function $\xdag$. In applications, one would typically expect
only polynomial decay of the Fourier coefficients of $\xdag-\unknown_0$ which corresponds to finite Sobolev smoothness
instead of infinite smoothness. Therefore, it is desirable to consider also functions $\Lambda$
which decay to $0$ more slowly than $t\mapsto t^{\nu}$. For exponentially
decaying singular values the logarithmic functions
\[
\Lambda(t) = (-\ln t)^{-p}
\]
with a parameter $p>0$ are a natural choice corresponding
to a polynomial decay of the Fourier coefficients of $\xdag-\unknown_0$. (Here we always assume that the operator is scaled such that 
$\|\linOp^*\linOp\|\leq \exp(-1)$ or alternatively use a dilated
version of the above function $\Lambda$.)
The importance of logarithmic source conditions for 
nonparametric instrumental regression is also pointed out in \cite{BluCheKri:07} and \cite{HorLee:07}.

\subsection{Variational source conditions for Banach spaces}\label{sec:source:var}

In our analysis we will not restrict ourselves to Hilbert spaces, but study the more general situation where $\mathcal{X}$
is a Banach space, which we assume in the following.
Note that in this case the operator $\linOp^*\linOp$ maps from $\mathcal{X}$ to the dual space $\mathcal{X}'$,
so even integer powers of $\linOp^*\linOp$ are not well-defined. Therefore, spectral source conditions as introduced above must be generalized. For this purpose we use variational methods which have been explored in regularization theory recently in a number 
of papers. We will prove convergence results with these methods in terms of the Bregman distance in $\mathcal{X}$
with respect to the convex functional $\mathcal{R}$.

Let $\xdag_*\in \partial \calR(\xdag)$ be a fixed element of the subdifferential of $\calR$ at $\xdag$ (i.e.\ $\xdag_*=\calR'[\xdag]$, if $\calR$ is differentiable at $\xdag$). Then the Bregman distance with respect to $\mathcal{R}$ and $\xdag_*$ is defined as
\begin{equation}\label{eq:defi_Bregman}
\Delta(\unknown,\xdag):= \calR(\unknown)-\calR(\xdag)-\lsp \xdag_*,\unknown-\xdag\rsp.
\end{equation}
Here $\lsp \cdot,\cdot\rsp$ denotes the classical dual paring $\lsp \calX', \calX \rsp$, i.e. $\lsp \xdag_*,\unknown-\xdag\rsp$
is the evaluation of the functional $\xdag_*$ at $\unknown-\xdag$. Hence, the Bregman distance measures how much
the linearization of $\calR$ at $\xdag$ and $\calR$ differ at the point $\varphi$. This is illustrated in Figure \ref{fig:bregman}. 
For strictly convex $\calR$ we have $\Delta(\unknown,\xdag) = 0$ if and only if $\unknown=\xdag$. The Bregman distance $\Delta$ is nonnegative and convex in the first argument, but it does not define
a metric since it is neither symmetric nor does it satisfy the triangle inequality in
general. However, Bregman distances provide a generalization of the simpler case, where $\mathcal{X}$ is a Hilbert space and
$\calR(\unknown)=\|\unknown-\unknown_0\|_{\mathcal{X}}^2$
for some $\unknown_0\in\mathcal{X}$. Because, in this situation
\begin{floatingfigure}[r]{7cm}
\includegraphics[width=7cm]{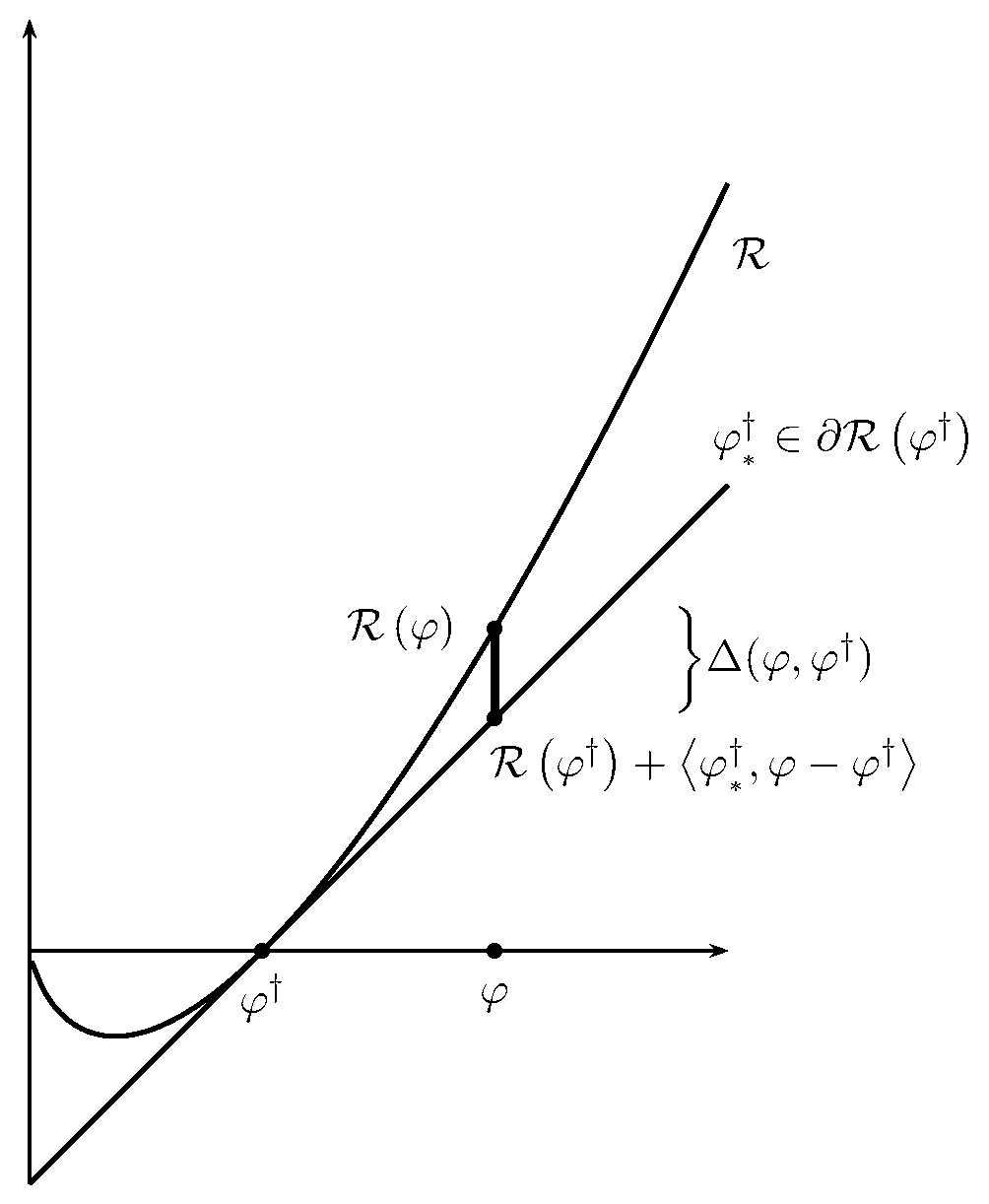}
\caption{Bregman distance}
\label{fig:bregman}
\vspace{5mm}
\end{floatingfigure}
\[
\Delta(\unknown,\xdag) = \|\unknown-\xdag\|_{\mathcal{X}}^2.
\]
Although, Bregman distances are in general not metrics they have meaningful interpretations in some Banach space settings.
If $\mathcal{X} = L^1(D)$ and $\mathcal{R}(\varphi) = \int_D \varphi(x) \ln(\varphi(x))\, dx $
(entropy regularization), then $\Delta(\unknown,\xdag)$ can be bounded from below by $\|\unknown-\xdag\|_{L^1}^2$ (see e.g.\linebreak
\cite{resmerita:05}), i.e.\ the error bounds formulated in the next theorem can be interpreted as 
bounds with respect to the squared $L^1$ norm. Our framework also allows the
incorporation of convex constraints by setting $\calR(\unknown):=\infty$ if
$\unknown$ does not belong to some convex set $\mathcal{C}$. Obviously, this does
not change $\Delta$ in $\mathcal{C}$.

Following \cite{KH:10} we formulate the source condition as a variational inequality
\begin{equation}\label{eq:sc_mult}
\lsp \xdag_*,\xdag-\unknown\rsp \leq \beta \Delta(\unknown,\xdag)^{1/2}
\Lambda\paren{\frac{\|\Op'[\xdag](\unknown-\xdag)\|^2}{\Delta(\unknown,\xdag)}}
\qquad \mbox{for all }\unknown\in \mathfrak{B}.
\end{equation}
Again, this is a generalization of the Hilbert space case. It is shown in \cite{KH:10} that if $\mathcal{X}$ is a Hilbert space,
$\mathcal{R}(\unknown)=\|\unknown-\unknown_0\|^2$ and $(\Lambda^2)^{-1}$ is convex, the classical source condition \eqref{eq:classical_sc} implies the variational one \eqref{eq:sc_mult}.

Let us close this section with a technical remark. Note that if $\mathfrak{B}$ is chosen such that $\xdag$ is on the boundary of $\mathfrak{B}$, then possibly $\Lambda$ can be chosen smaller than in the case where $\xdag$ is in the interior of $\mathfrak{B}$. Theorem \ref{theo:main} yields that this may lead to faster rates of convergence. Hence, a convex constraint on the regression function can improve estimation. To captures this fact it is important that, opposed to the formulation in \cite{KH:10}, no absolute values appear on the left hand side of \eqref{eq:sc_mult}. A typical example where $\xdag$ is on the boundary of $\mathfrak{B}$ is the assumption that $\xdag$ is a positive function.

\section{Convergence results}\label{sec:convergence}
Let $\calX$ be a Banach space, $\calY$ a Hilbert space, $\mathfrak{B}\subset \calX$ convex and $\xdag\in\mathfrak{B}$ a root of the operator $\Op:\mathfrak{B}\to \calY$:
\begin{equation}\label{eq:defi_xdag}
\Op(\xdag)=0.
\end{equation}
Assume that $\Op$ is approximated by a series of estimators
\[
\Fhat_n:\mathfrak{B}\to \Yhat_n
\]
which maps to some (possibly finite-dimensional and/or data dependent)
Hilbert space  $\Yhat_n$. $\Op$ and all $\Fhat_n$ are assumed to be Gateaux
differentiable on $\mathfrak{B}$ with linear derivatives
$\Op'[\varphi]$ and $\Fhat_n'[\varphi]$, which are ``bounded with respect
to $\Delta$'' in the sense that
\begin{equation}\label{eq:bounded_by_delta}
\sup_{\{\tilde{\varphi}\in\mathfrak{B}:\Delta(\tilde{\varphi},\varphi)\neq 0\}}
\|\Op'[\varphi](\tilde{\varphi}-\varphi)\|^2/\Delta(\tilde{\varphi},\varphi)<\infty \quad \text{and} \quad \Op'[\varphi](\tilde{\varphi}-\varphi)\neq 0
\end{equation}
whenever $\Delta(\tilde{\varphi},\varphi)\neq 0$ and analogously for all $\Fhat_n$.
Now we can state the main theorem of this paper, which is proved in Appendix \ref{sec:proof}:
\begin{theorem}\label{theo:main}
Let \eqref{eq:sc_mult} hold true with a concave $\Lambda$ for which $t\mapsto \sqrt{t}/\Lambda(t)$ is monotonically increasing. Assume that the sequence $\Fhat_n$ has the following convergence properties:
\begin{subequations}
\begin{eqnarray}
&&\|\Fhat_n(\xdag)\| = O_p(\delta_{n}),\\
&&\paren{\left|\sup_{\unknown\in \mathfrak{B}} \frac{\|\Op'[\xdag](\unknown-\xdag)\|^2-\|\Fhat_n'[\xdag](\unknown-\xdag)\|^2}
{\Delta(\unknown,\xdag)}\right|}^{1/2} = O_p(\gamma_{n}),\\
&&\begin{split}\label{eq:probtc}
&P\{\|\Fhat_n(\unknown_1)-\Fhat_n(\unknown_2)-\Fhat_n'[\unknown_2](\unknown_1-\unknown_2)\| > \eta \|\Fhat_n(\unknown_1)-\Fhat_n(\unknown_2)\|\\
&\hspace{78mm} \text{ for some } \unknown_1,\unknown_2 \in \mathfrak{B} \} \rightarrow 0.
\end{split}
\end{eqnarray}
\end{subequations}
Here $\eta$ must be sufficiently small, such that $4 \eta (1 + \eta)(1 - \eta)^{-3} < q^{-3/2}$.
Suppose that the convex minimization problems
\eqref{eq:defi_method} are uniquely solvable for every $\Fhat_n$ (see Remark \ref{rem:ex_unique}
for sufficient conditions), i.e.\ the method is well defined. Further assume that $\alpha_0>\max(\Theta^{-1}(\delta_n),\gamma_n^2)$ and that $\alpha_k \leq q \alpha_{k+1}$ for all $k$ with a constant $q > 1$. Let the iteration be stopped at the smallest index $\capit_n \in \mathbb{N}_0$ for which
\begin{equation}\label{eq:apriori_stopping}
\alpha_{\capit_n +1} \leq \max(\Theta^{-1}(\delta_{n}),\gamma_{n}^2)
\,,\qquad \mbox{where } \Theta(t):=\sqrt{t}\Lambda(t).
\end{equation}
Then
\[
\Delta(\widehat{\unknown}_{\capit_n},\xdag) = O_p\left(\Lambda^2\left(\max(\Theta^{-1}(\delta_{n}),\gamma_{n}^2)\right)\right)
\]
\end{theorem}

\noindent \textbf{Remarks:}
\begin{enumerate}
\item\label{rem:ex_unique} Sufficient conditions for uniqueness of solutions to the minimization problem
\eqref{eq:defi_method} are strict convexity of $\mathcal{R}$
or injectivity of $\Fhat_n'[\hatxkm]$.

\item Sufficient conditions for existence are reflexivity of $\mathcal{X}$,
weak closedness of $\mathfrak{B}$, and the boundedness of the sets
$\{\unknown\in\mathfrak{B}:\mathcal{R}(\unknown)\leq R\}$ in $\mathcal{X}$
for any $R\in\mathbb{R}$. This is a standard argument: If
$(\unknown_n)$ is a minimizing sequence, it must be bounded due to our last
condition. Since $\mathcal{X}$ is reflexive, there exists a weakly
convergent subsequence, and by weak
closedness of $\mathfrak{B}$ a weak limit point $\unknown_*\in\mathfrak{B}$.
Since the Tikhonov functional is convex and lower semi-continuous, it
is also weakly lower semi-continuous, and hence $\unknown_*$ is a minimizer.

\item Note that if $\mathcal{X}$ is a Hilbert space and $\Fhat_n$ Fr\'echet differentiable, then $\|\Op'[\xdag](\unknown-\xdag)\|^2-\|\Fhat_n'[\xdag](\unknown-\xdag)\|^2
\leq \|\Op'[\xdag]^*\Op'[\xdag] - \Fhat_n'[\xdag]^*\Fhat_n'[\xdag]\|\; \|\unknown-\xdag\|^2$,
so  $\gamma_n\leq \|\Op'[\xdag]^*\Op'[\xdag]-\Fhat_n'[\xdag]^*\Fhat_n'[\xdag]\|^{1/2}$.

\item The bound on the Taylor remainder of $\Fhat_n$
\begin{equation}\label{eq:tcn}
\|\Fhat_n(x)-\Fhat_n(y)-\Fhat_n'[y](x-y)\|\leq \eta \|\Fhat_n(x)-\Fhat_n(y)\|,
\end{equation}
used in \eqref{eq:probtc} is known as the tangential cone condition.
This condition is commonly used in the analysis of regularization methods
for nonlinear ill-posed problems, see \cite{KNS:08}. 
The right hand side of \eqref{eq:tcn} may be replaced by $\|F'[y](x-y)\|$
(see \eqref{eq:coro_tc} below), and in this form it corresponds to 
Assumption 2 in \cite{CCLN:11}.
\end{enumerate}

\begin{corollary} \label{cor:rates}
Let the assumptions of Theorem \ref{theo:main} hold true.
\begin{enumerate}
\item If $\Lambda(t) = t^{\mu}$ for
some $\mu\in (0,1/2]$  (H\"older-type source conditions), then
\begin{equation}\label{eq:HoelderRate}
\Delta(\widehat{\unknown}_{\capit},\xdag)^{1/2}
=O_p\left(\max(\delta_n^{2\mu/(2\mu+1)}, \gamma_n^{2\mu})\right)\,.
\end{equation}
\item  If $\Op$ is scaled such that
$\|\Op'[\xdag](\unknown-\xdag)\|^2/\Delta(\unknown,\xdag)\leq \frac{1}{2}$
and $\Lambda(t) = (- \ln t)^{-p}$ for some $p>0$ (logarithmic source conditions), then
\begin{equation}\label{eq:logarithmicRate}
\Delta(\widehat{\unknown}_{\capit},\xdag)^{1/2}
= O_p\left((-\ln \max(\delta_n,\gamma_n))^{-p}\right)
\end{equation}
for all $\delta_n,\gamma_n$ sufficiently small.
\end{enumerate}
\end{corollary}

Let us discuss some properties of the method. First of all it is a local method
like any Newton method. 
Convergence is only guaranteed if the initial guess $\varphi_0$ is sufficiently close to the true solution $\xdag$. 
How close it has to be depends on the special problem, i.e. the operator $\Op$.
This property appears in the assumptions 
\eqref{eq:sc_mult} and \eqref{eq:probtc} in Theorem \ref{theo:main}.

We emphasized in the introduction that the method requires only solutions of convex minimization problems. 
Therefore, it does not suffer from the problem of multiple local minima 
which frequently occur in nonlinear Tikhonov regularization \eqref{eq:nltikh} 
and make it hard to find the actual minimum.

Unlike for nonlinear Tikhonov regularization our theoretical results do not require the strong assumption that we can always find
the minimum of a functional with an arbitrary number of local minima. In turn we have to assume
\eqref{eq:probtc}, which is usually hard to check. Although, rigorous proofs for \eqref{eq:probtc} are often 
missing, it seems to hold in many cases at least in a neighborhood of $\xdag$.

An important advantage for the numerical implementation is that a lot of efficient algorithms converging always towards 
the true solution are known for convex minimization problems. 
The error of these minimization algorithms plays a minor role compared to the regularization error for the applications of Section \ref{sec:example}. We refer to \cite{HL07} for a detailed discussion of the interplay of these errors in other applications.

\section{Examples revisited}\label{sec:revisited}

The assumptions of Theorem \ref{theo:main} and Corollary \ref{cor:rates} are rather abstract and need some explanations concerning the application 
 to the nonparametric regression with independent instrument \eqref{eq:instr_model}. They are applicable 
in a similar way to the nonparametric quantile regression \eqref{eq:quantile_regression}. In \eqref{eq:nonlininstreg:opeq} the
operator $\Op$ for the regression with independent instrument is an integral operator with a kernel composed by marginals of 
$f_{YZW}$. Hence, an estimator $\widehat f_{YZW}$ yields an estimator of the kernel and thereby of $\Ophat$.

Condition \eqref{eq:bounded_by_delta} that all $\Ophat_n'[\unknown]$ must be bounded with respect to the Bregman distance is fulfilled if the derivatives of $\Op$ are bounded according to \eqref{eq:bounded_by_delta}, the estimation of $f_{YZW}$ is strongly consistent and $n$ is large enough. Strong consistency is established for many density estimators. The boundedness of $\Op$ with respect to the Bregman distance is reasonable. It holds if the partial derivative of the joint density $\frac{\partial}{\partial y}f_{YZW}$ is bounded for the operator \eqref{eq:nonlininstreg:opeq} or if $f_{YZ|W}$ is bounded for the operator \eqref{eq:nonlininstreg:opeq2} and the Bregman distance is bounded from below by the power of a norm. As mentioned in Section \ref{sec:source:var}, the latter is for example the case for quadratic and maximum entropy penalty.

It can be argued with strong consistency as well that the probabilistic tangential cone condition \eqref{eq:probtc} holds if the exact operator $\Op$ fulfills the tangential cone condition \eqref{eq:tcn}. But, it is known that a verification, whether or not the tangential cone condition is true, is often difficult for a given operator.

In analogy to \eqref{eq:nonlininstreg:opeq2} the operator
\[
\tilde\Op(\varphi)(u,w) := \left(\begin{array}{c}
 \Prob \big(Y - \varphi(Z) \leq u \big) - \Prob \big(Y - \varphi(Z) \leq u | W = w \big)\\
\Ex[\varphi(Z) - Y]
\end{array}\right)
\]
can be considered for model \eqref{eq:instr_model}. With this operator the conditions \eqref{eq:bounded_by_delta} and \eqref{eq:tcn} are more explicitly assumptions on the primitives of the model.

In the rates for H\"older source conditions \eqref{eq:HoelderRate} in Corollary \ref{cor:rates}, $\delta$ has a smaller exponent than $\gamma$. However $\delta$ does not necessarily dominate the convergence. In the nonparametric instrumental regression $\delta$ corresponds to the estimation of a density, while $\gamma$ is determined by the estimation of a partial derivative of that density. Hence, $\gamma$ decays usually slower than $\delta$. Which of the terms $\delta^{2\mu/(2\mu+1)}$ or $\gamma^{2\mu}$ dominates the convergence depends on the properties of the special problem, namely the number of instruments and covariates as well as on the smoothness of the density and the initial error $\xdag-\unknown_0$.

The situation becomes clearer in the case of logarithmic source conditions. If the kernel of the operator is analytic, but the initial error in the regression function is not smooth or has only finite H\"older smoothness, merely a logarithmic rate of convergence can be expected. As discussed in Section \ref{sec:source:impact} this situation can occur in many applications. Even for estimating an analytic density a nonparametric density estimator will attain only a polynomial rate in $n$. Due to \eqref{eq:logarithmicRate} our estimator for $\varphi$ will end up asymptotically with the logarithmic rate $(-\ln(n))^{-p}$.

\section{Numerical simulations}\label{sec:simulations}

In this section we  present some numerical simulations for nonparametric instrumental regression with independent binary instrument and real-valued continuous explanatory and dependent variables. This leads to the nonlinear operator equation \eqref{eq:nonlinbininstreg:opeq}. Our simulations show that the solution computed by the method \eqref{eq:defi_method} approximates the exact solution. As mentioned above, due to dimensionality, the regression function cannot be identified with a binary instrument if the standard regression model \eqref{eq:standard_instr} is used.

In our simulations we choose  $Y$ as real valued, $Z $ with values in $ [0,1]$ and $W $ with values in $ \{0,1\}$. We assume the regression function is 
\[
 \varphi^\dag(z) = \frac{1}{6}\sin(2 \pi (z+0.25))+0.41\,,\qquad z\in [0,1]\,.
\]
Moreover, we take $w_0 = P(W=0) = 2/3$ and $w_1 = P(W=1) = 1/3$. To make $Z$ endogenous, let us choose the error term as $(U|Z=z, W=w) \sim \mathcal{N}(\mu_w(z), 0.09^2)$ with 
$\mu_0(z) := 0.2 z - 0.1$ and $\mu_1(z) := 0.25 z - 0.125$.
The functions $\mu_0(z)$ and $\mu_1(z)$ describe precisely the correlation between the explanatory variable and the error term, which should be removed using the information contained in the instrumental variable. 
Although $U$ varies with $Z$ and $W$ the condition $W \independent U$ can be assured by a proper choice of $f_{Z,W}(z,w)$.  We write the joint density as
\begin{equation}\label{eq:fyzw1}
\begin{split}
f_{YZW}(y,z,0) = f_{ZW}(z,0)\frac{1}{0.09\sqrt{2\pi}}\exp\left(-\frac{1}{2}\left(\frac{y-\varphi^{\dag}(z)-\mu_0(z)}{0.09}\right)^2\right),\\
 f_{YZW}(y,z,1) = f_{ZW}(z,1)\frac{1}{0.09\sqrt{2\pi}} \exp\left(-\frac{1}{2}\left(\frac{y-\varphi^{\dag}(z)-\mu_1(z)}{0.09}\right)^2\right).
\end{split}
\end{equation}
Now $f_{ZW}$ has to be determined such that $W$ and $U$ are independent, which is equivalent to \eqref{nullstelle}. Let us show that setting $f_{ZW}(z,1) := 0.625 f_{ZW}(1.25 z - 0.125,0)$ achieves this. With a substitution of variables we compute
\begin{equation*}
\begin{split}
& \int w_1 f_{YZW}(u + \varphi^{\dag}(z), z, 0) \,dz\\
= & \int \frac{1}{3} f_{ZW}(z,0)\frac{1}{0.09\sqrt{2\pi}}\exp\left(-\frac{1}{2}\left(\frac{u-0.2z+0.1)}{0.09}\right)^2\right)\,dz\\
= & \int \frac{1.25}{3} f_{ZW}(1.25v - 0.125,0)\frac{1}{0.09\sqrt{2\pi}}\exp\left(-\frac{1}{2}\left(\frac{u-0.25v+0.125)}{0.09}\right)^2\right)\,dv\\
= & \int w_0 f_{YZW}(u + \varphi^{\dag}(v), v, 1)\,dv.
\end{split}
\end{equation*}
This shows that \eqref{nullstelle} holds with our definition of $f_{ZW}(z,1)$ what ever $f_{ZW}(z,0)$ looks like. Here we take it to  be normally distributed with variance $0.3^2$ and expectation $1/2$ truncated to the interval $[0,1]$, i.e.
\[
 f_{ZW}(z,0) := a\exp\left(-\frac{1}{2}\left(\frac{z-1/2}{0.3}\right)^2\right),
 \qquad z\in [0,1]
\]
with some scaling factor $a$ chosen such that $\int_0^1 f_{ZW}(z,0) dz = 2/3$. 
By this construction, the error term also meets the condition 
$\mathbb{E}U = 0$ of the regression model \eqref{eq:instr_model}: 
To see this, note that 
$f_{ZW}(\cdot,0)$ and $f_{ZW}(\cdot,1)$ are even, 
while $\mu_0$ and $\mu_1$ are odd functions with respect to the point $0.5$. 
Hence, 
\begin{equation*}
\begin{split}
\mathbb{E}U &= \int w_0 f_{Z,W}(z,0)\mathbb{E}(U|Z=z, W=0) + w_1 f_{Z,W}(z,1)\mathbb{E}(U|Z=z, W=1) \,dz\\
&= \int w_0 f_{Z,W}(z,0)\mu_0(z) + w_1 f_{Z,W}(z,1)\mu_1(z)\,dz = 0.
\end{split}
\end{equation*}

This construction allows an easy formulation of how the solution of a nonparametric regression without instrumental variable and without noise would look like: $\widetilde\varphi(z) = w_0\mu_0(z) + w_1\mu_1(z) + \varphi^\dag$
\begin{center}
\includegraphics[width=0.5\textwidth]{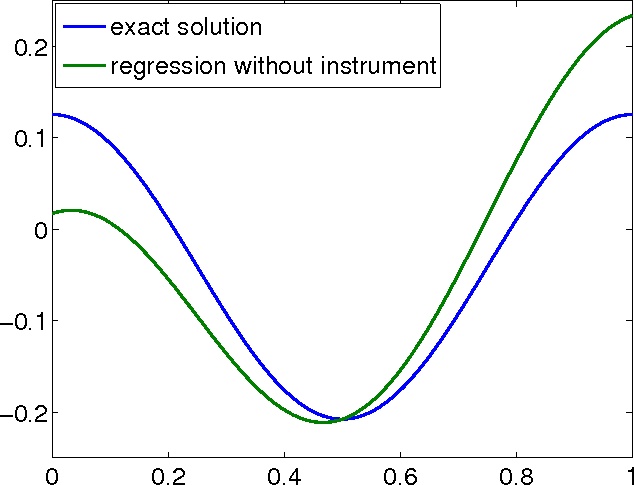}
\captionof{figure}{Necessity of the instrument: A standard regression would  asymptotically yield a curve which is far away from the true solution $\varphi^\dag$.}
\end{center}

To approximately solve the integral operator equation \eqref{eq:nonlinbininstreg:opeq} by the method \eqref{eq:defi_method}
we discretized the domain $[0,1]\times[0,1]\times\{0,1\}$ by
$256 \times 256 \times 2$ points and chose the regularization parameters by
$\alpha_0 = 1$ and $\alpha_{n+1} = 0.9 \alpha_n$. The iteration was stopped using Lepski\u{\i}'s principle as in Bauer, Hohage \& Munk \cite{BauHohMun:09}.
The initial guess was chosen as the constant function $E[Y]$. For a
first test we used the exact density $f_{YZW}$, which actually has to be
estimated from the data, of course.
The $L^2$-error was reduced from $0.1294$ to $0.0028$. The remaining error is
due to discretization noise. This suggests that the example is identifiable
and can be solved by the method \eqref{eq:defi_method}. Compared to the error for 
densities estimated from simulated data below, the observed discretization error is very small. Hence, the 
discretization is fine enough and discretization error is insignificant for our simulations. 
The singular values of $\Op'[\varphi^\dag]$ are shown in Figure \ref{fig:svd}.
They exhibit an exponential decay, so according to  Corollary \ref{cor:rates}
we can only expect slow rates of convergence.\\[5mm]
\begin{minipage}{0.49\textwidth}
 \includegraphics[width=\textwidth]{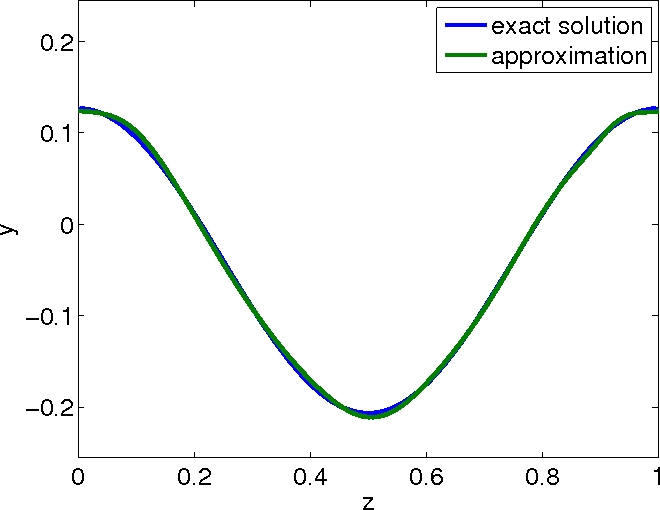}
 \captionof{figure}{Result of the iterative inversion 
 using the exact density $f_{YZW}$.}
\end{minipage}\quad
\begin{minipage}{0.49\textwidth}
  \includegraphics[width=\textwidth]{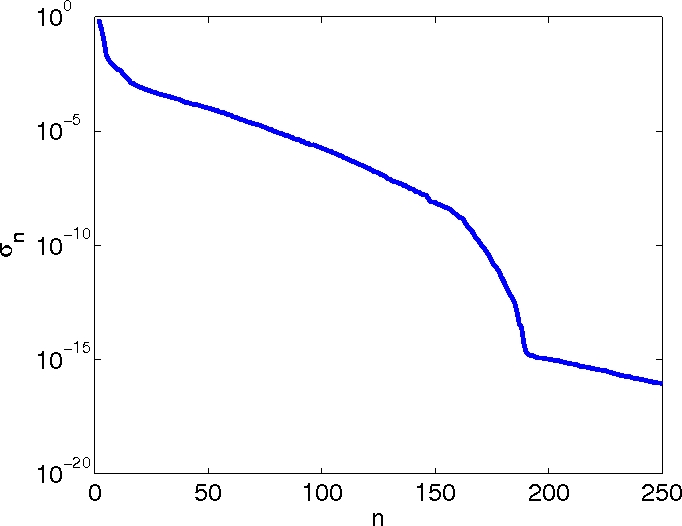}
   \captionof{figure}{\label{fig:svd}Singular values of $\Op'[\varphi^\dag]$}
 \end{minipage}\\[5mm]
In further tests the algorithm was evaluated for finite samples of $(Y,Z,W)$ with $10^3$, $10^4$ and $10^5$ points. Given such a sample, the joint density $f_{YZW}$ was estimated non-parametrically by the kernel density estimator developed by Botev, Grotowski \& Kroese \cite{BotGroKro:10}. Afterwards again \eqref{eq:nonlinbininstreg:opeq} was solved, but the exact density was replaced by the estimated one. We made $1000$ samples for each tested sample size.
The following table and histograms in Fig.~\ref{fig:hist3}--\ref{fig:hist5}
show the $L^2$-errors of the approximate
solution normed by the error of the initial guess (i.e. the error of the initial guess is $1$). It can be seen that small samples produce unwanted outliers, but the method becomes reliable 
when the sample size is large enough. Fig.~\ref{fig:rec3}--\ref{fig:rec5}
show median reconstructions for each sample size.
The results demonstrate that our method computes an asymptotically correct
estimator of the regression function $\varphi^{\dag}$ with an endogeneous
explanatory variable $Z$ using only a binary instrument $W$. 

the exact solution is $0$ and the error of the initial guess is $1$. It can be seen that small samples produce unwanted outliers, but that the method becomes reliable, when the sample is large enough.

\begin{center}
\begin{tabular}{|l||l|ll|l|l|l|}
\hline
sample size $N$ & mean & quantiles & $p = 0.25$ & $p = 0.5$ & $p = 0.75$ & $p = 0.9$\\ \hline
$10^3$ & 0.6159 & & 0.4057 & 0.5751 & 0.7921 & 0.9575\\ \hline
$10^4$  & 0.3694 & & 0.2496 & 0.3524 & 0.4574 & 0.5729\\ \hline
$10^5$ & 0.3264 & & 0.2592 & 0.3278 & 0.3882 & 0.4610\\ \hline
\end{tabular}
\end{center}
\begin{minipage}{0.49\textwidth}
 \includegraphics[width=\textwidth]{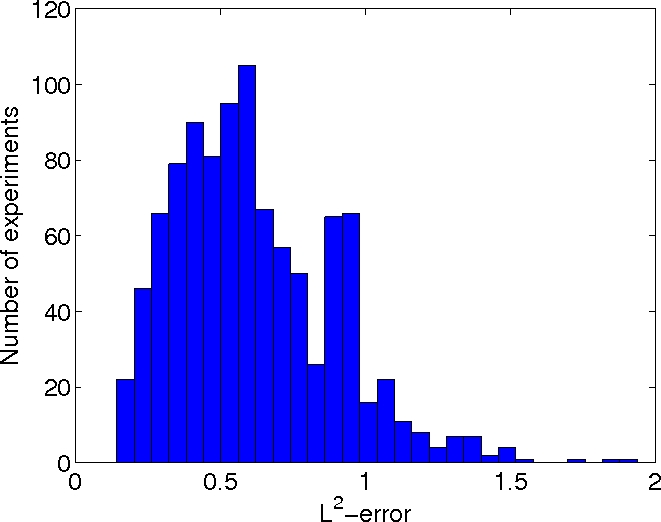}
 \captionof{figure}{\label{fig:hist3}$L^2$ error for sample size $N=10^3$}
\end{minipage}\quad
\begin{minipage}{0.49\textwidth}
 \includegraphics[width=\textwidth]{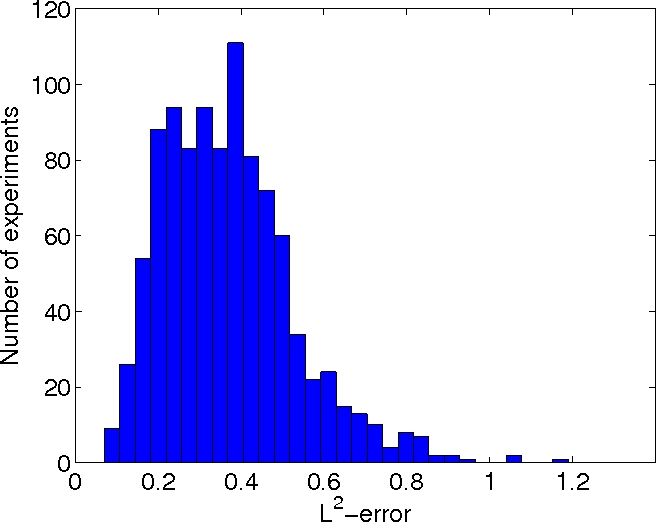}
 \captionof{figure}{\label{fig:hist4}$L^2$ error for sample size $N=10^4$}
\end{minipage}\\
\begin{minipage}{0.49\textwidth}
\includegraphics[width=\textwidth]{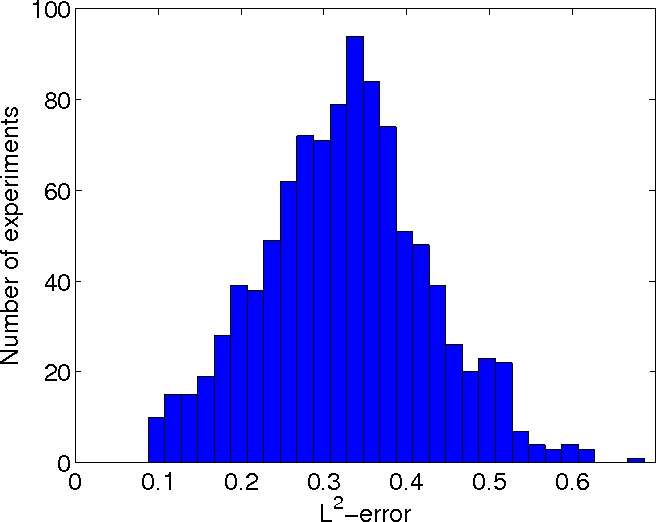}
\captionof{figure}{\label{fig:hist5}$L^2$ error for sample size $N=10^5$}
\end{minipage}\quad
\begin{minipage}{0.49\textwidth}
 \includegraphics[width=\textwidth]{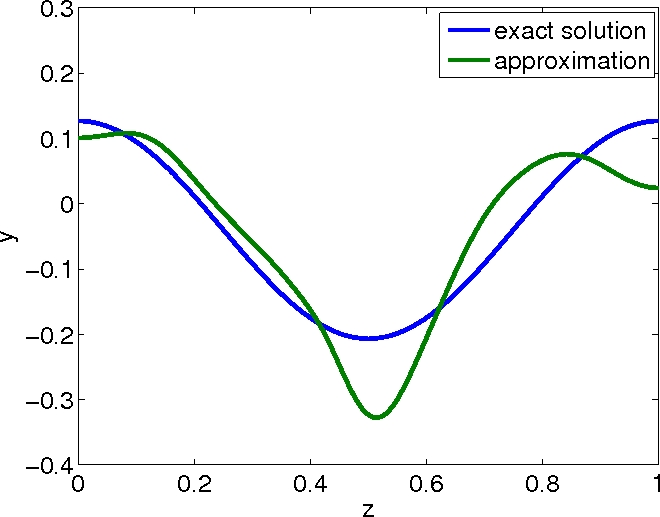}
 \captionof{figure}{\label{fig:rec3}Median reconstruction, $N=10^3$ }
\end{minipage}\\
\begin{minipage}{0.49\textwidth}
\includegraphics[width=\textwidth]{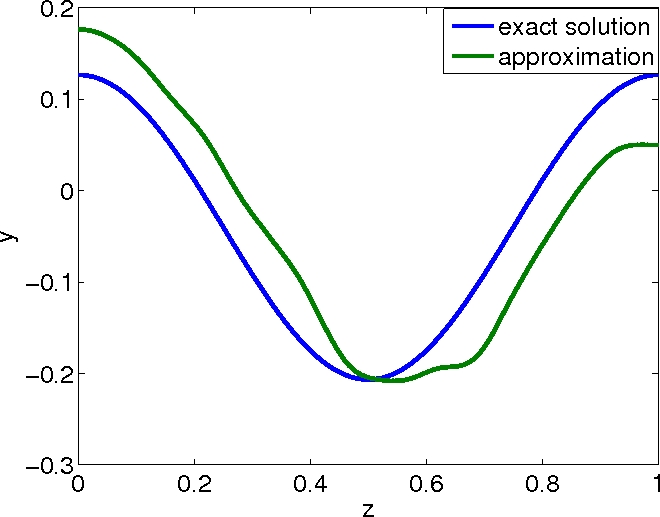}
\captionof{figure}{\label{fig:rec4}Median reconstruction, $N=10^4$}
\end{minipage}\quad
\begin{minipage}{0.49\textwidth}
 \includegraphics[width=\textwidth]{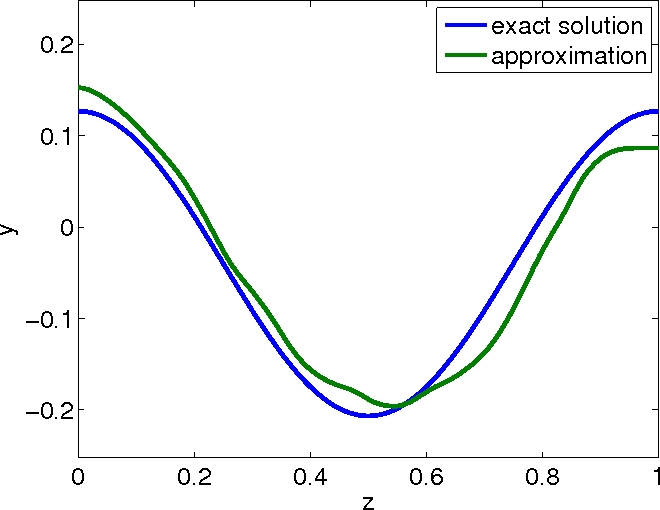}
 \captionof{figure}{\label{fig:rec5}Median reconstruction, $N=10^5$}
\end{minipage}\\[5mm]

\appendix
\section{Proof of the main theorem}\label{sec:proof}
Before we come to the proof of Theorem \ref{theo:main}, let us first formulate a result with deterministic error in the operator. We assume that $\Op$ is approximated by some deterministic operator
\[
\Fhat:\mathfrak{B}\to \Yhat.
\]
Let both $\Op$ and $\Fhat$ be Gateaux
differentiable on $\mathfrak{B}$ with derivatives
$\Op'[\varphi]$ and $\Fhat'[\varphi]$, which are ``bounded with respect
to $\Delta$'' in the sense that
$\sup_{\{\tilde{\varphi}\in\mathfrak{B}:\Delta(\tilde{\varphi},\varphi)\neq 0\}}
\|\Op'[\varphi](\tilde{\varphi}-\varphi)\|^2/\Delta(\tilde{\varphi},\varphi)<\infty$ and $\Op'[\varphi](\tilde{\varphi}-\varphi)\neq 0$ whenever
$\Delta(\tilde{\varphi},\varphi)\neq 0$ and analogously for $\Fhat$.
The error of the approximation is described by:
\begin{subequations}\label{eqs:defi_delta}
\begin{eqnarray}
\label{eq:defi_delta1}
\delta&:=& \|\Fhat(\xdag)\|,\\
\label{eq:defi_delta2}
\gamma&:=& \paren{\left|\sup_{\{\unknown\in \mathfrak{B}:\Delta(\unknown,\xdag)\neq 0\}}
\frac{\|\Op'[\xdag](\unknown-\xdag)\|^2-\|\Fhat'[\xdag](\unknown-\xdag)\|^2}
{\Delta(\unknown,\xdag)}\right|}^{1/2}.
\end{eqnarray}
\end{subequations}
Moreover, we assume that the tangential cone condition
\begin{equation}\label{eq:tc}
\|\Fhat(x)-\Fhat(y)-\Fhat'[y](x-y)\|\leq \eta \|\Fhat(x)-\Fhat(y)\|,
\end{equation}
holds for all $x,y$ in some neighborhood of $\mathcal{B}$.

\begin{lemma}\label{theo:main_lemma}
Assume that \eqref{eq:sc_mult}, \eqref{eqs:defi_delta} and \eqref{eq:tc} hold true
with $\eta$ sufficiently small, such that
\begin{align}\label{eq:small_eta}
4 \eta (1 + \eta)(1 - \eta)^{-3} < q^{-3/2}.
\end{align}
Further assume that the convex minimization problems
\eqref{eq:defi_method} are uniquely solvable and
 that the iteration is stopped at the smallest
index $\capit\in\mathbb{N}_0$ for which
\begin{equation}\label{eq:apriori_stopping2}
\alpha_{\capit+1} \leq \max(\Theta^{-1}(\delta),\gamma^2)
\,,\qquad \mbox{where } \Theta(t):=\sqrt{t}\Lambda(t).
\end{equation}
In addition it should hold that $\alpha_0>\max(\Theta^{-1}(\delta),\gamma^2)$ and $\alpha_k \leq q \alpha_{k+1}$ 
for all $k$ with a constant $q > 1$. Moreover, let
$\Lambda$ be concave and assume that $t\mapsto \sqrt{t}/\Lambda(t)$ is monotonically increasing.

Then there exists a constant $C>0$ independent of
the $\Fhat$ such that
\begin{eqnarray}\label{eq:estim_Delfin}
\Delta(\widehat{\unknown}_{\capit},\xdag) &\leq&
C \paren{\Lambda\left(\max(\Theta^{-1}(\delta),\gamma^2)\right)}^2\,.
\end{eqnarray}
\end{lemma}

\proof
Let us introduce the following notation:
\begin{eqnarray*}
 &&\T:= \Op'[\xdag],\qquad
 \Tdag:= \Fhat'[\xdag], \qquad
 \Tkm:=\Fhat'[\hatxkm],\\
 &&\Delk:=\Delta(\hatxk,\xdag),\qquad  \ek:=\hatxk-\xdag.
 \end{eqnarray*}

From the optimality condition \eqref{eq:defi_method} with $\unknown=\xdag$  we find that
\begin{equation}\label{eq:optimality_cond}
\begin{split}
&\|\Tkm(\hatxk-\hatxkm) + \Fhat(\hatxkm)\|^2 + \alk \calR(\hatxk)\\
&\qquad \leq \|\Tkm(\xdag-\hatxkm) + \Fhat(\hatxkm)\|^2 + \alk \calR(\xdag).
\end{split}
\end{equation}
From the definition \eqref{eq:defi_Bregman} of the Bregman distance
and the source condition \eqref{eq:sc_mult} we obtain
\begin{equation}\label{eq:use_of_sc_mult}
\calR(\xdag)-\calR(\hatxk) = \lsp \xdag_*,\xdag-\hatxk\rsp - \Delk
\;\leq\; \beta \Delk^{1/2}\Lambda\paren{\frac{\|\T \ek\|^2}{\Delk}} - \Delk.
\end{equation}
Plugging this into \eqref{eq:optimality_cond} yields
\begin{equation}\label{eq2:aux1}
\begin{split}
&\|\Tkm(\hatxk-\hatxkm) + \Fhat(\hatxkm)\|^2
+ \alk \Delk \\
&\qquad \leq \|\Tkm(\xdag-\hatxkm) + \Fhat(\hatxkm)\|^2
+ \beta \alk \Delk^{1/2}\Lambda\paren{\frac{\|\T \ek\|^2}{\Delk}} \,.
\end{split}
\end{equation}
Note that the tangential cone condition \eqref{eq:tc} implies
\begin{equation}\label{eq:coro_tc}
(1+\eta)^{-1}\|\Tdag\ek\| \leq \|\Fhat(\hatxk)-\Fhat(\xdag)\|\leq (1-\eta)^{-1} \|\Tdag \ek\|\,.
\end{equation}
To estimate the first term on the left hand side of \eqref{eq2:aux1}
we use \eqref{eq:tc} and \eqref{eq:coro_tc} to get that
\begin{eqnarray*}
\lefteqn{\|\Fhat(\hatxk)\| - \|\Tkm(\hatxk-\hatxkm) + \Fhat(\hatxkm)\|}\\
&\leq& \|\Tkm(\hatxk-\hatxkm) + \Fhat(\hatxkm) - \Fhat(\hatxk)\| \\
&\leq& \eta \|\Fhat(\hatxkm) - \Fhat(\hatxk)\|\\
&\leq& \eta\|\Fhat(\hatxkm) - \Fhat(\xdag)\|+\eta\|\Fhat(\hatxk) - \Fhat(\xdag)\|\\
&\leq& \frac{\eta}{1-\eta} (\|\Tdag \ek\|+\|\Tdag \ekm\|).
\end{eqnarray*}
Together with $\||\Fhat(\hatxk)\| \geq \|\Fhat(\hatxk)-\Fhat(\xdag)\| -\delta
\geq (1+\eta)^{-1}\|\Tdag \ek\|-\delta$
 this yields
\begin{eqnarray*}
\|\Tkm(\hatxk-\hatxkm) + \Fhat(\hatxkm)\|
\geq \frac{(1-\eta)^2}{1-\eta^2}\|\Tdag\ek\| - \frac{\eta}{1-\eta} \|\Tdag\ekm\|
   - \delta.
\end{eqnarray*}
For the right hand side of \eqref{eq2:aux1} we get from \eqref{eqs:defi_delta}
and another application of \eqref{eq:tc} that
\begin{equation*}
\|\Tkm(\xdag-\hatxkm) + \Fhat(\hatxkm)\| \leq \eta\|\Fhat(\hatxkm)-\Fhat(\xdag)\| +\delta
\leq \frac{\eta}{1-\eta}\|\Tdag \ekm\| +\delta.
\end{equation*}
Plugging the last two inequalities into \eqref{eq2:aux1} and using the simple inequalities
$(a-b)^2\geq \frac{1}{2}a^2-b^2$ and $(a+b)^2\leq 2a^2+2b^2$ we obtain that
\begin{equation*}
\begin{split}
&\underbrace{\frac{1}{2}\left(\frac{(1-\eta)^2}{1-\eta^2}\right)^2}_{=:C_{\eta}}\left\| \Tdag\ek\right\|^2 + \alk  \Delk
 \leq \underbrace{\frac{4\eta^2}{(1-\eta)^2}}_{=:c_\eta}\left\|\Tdag\ekm\right\|^2
+ 4\delta^2
+  \beta \alk \Delk^{1/2}\Lambda\paren{\frac{\|\T \ek\|^2}{\Delk}}\,.
\end{split}
\end{equation*}
Using \eqref{eq:defi_delta2} and the monotonicity of $\Lambda$ we find that
$\Lambda\paren{\frac{\|\T \ek\|^2}{\Delk}}\leq
\Lambda\paren{\frac{\|\Tdag \ek\|^2}{\Delk}+\gamma^2}$.
Together with the stopping rule \eqref{eq:apriori_stopping2} this implies
\begin{equation}\label{eq:induction:inequality}
C_\eta \| \Tdag\ek\|^2 + \alk\Delk
\leq c_{\eta} \| \Tdag e_{\its-1}\|^2
+ 4\Theta(\alk)^2 + \beta \alk\Delk^{1/2}
\Lambda\paren{\frac{\| \Tdag\ek\|^2}{\Delk} + \alk}\,.
\end{equation}
We will show the following error bounds
\begin{subequations}
\begin{eqnarray}
\label{eq:estimt}
\| \Tdag\ek\|^2 &\leq& C_1\Theta(\alk)^2,\\
\label{eq:estimDelta}
\Delta(\hatxk,\xdag) &\leq& C_2\Lambda(\alk)^2
\end{eqnarray}
\end{subequations}
with
\begin{eqnarray*}
C_1&:=&\max\paren{\dfrac{\|\Tdag e_0\|^2}{\Theta(\alpha_0)^2},\, \dfrac{8}{C_\eta-2q^3c_\eta},\,
\dfrac{16\beta^2}{C_\eta+1},\, \dfrac{16\beta^2}{C_\eta^2}},\\
C_2&:=&\max\paren{\dfrac{\Delta(\varphi_0,\xdag)}{\Lambda(\alpha_0)^2},\,
2C_1c_\eta q^3+8,\, 16\beta^2,\, \frac{ 16\beta^2}{C_\eta}}.
\end{eqnarray*}
We will prove these claims by induction in $\its\leq \capit$. For $\its=0$ this
is arranged by the definitions of $C_1$ and $C_2$. For the induction step we distinguish two cases:\bigskip\\
\emph{Case 1:} $c_{\eta} \| \Tdag e_{\its-1}\|^2 + 4\Theta(\alk)^2
\geq \beta \alk\Delk^{1/2} \Lambda\paren{\frac{\| \Tdag\ek\|^2}{\Delk} + \alk}$.\\
Now by using the induction hypothesis \eqref{eq:estimt} equation \eqref{eq:induction:inequality} simplifies to
\[
C_\eta \| \Tdag\ek\|^2 + \alk\Delk \leq 2c_{\eta} C_1 \Theta(\alpha_{k-1})^2 + 8\Theta(\alk)^2 .
\]
We have $\Theta(\alpha_{k-1}) = (\alpha_{k-1})^{1/2} \Lambda(\alpha_{k-1}) \leq (q\alk)^{1/2} \Lambda(q \alk)$ as $\Lambda$ is monotonically increasing. While $\Lambda$ is concave and $\Lambda(0) = 0$ the definition of concavity implies $t \Lambda(x) \leq \Lambda(tx)$ for $0 \geq t \geq 1$. Now taking $x=q\alk$ and $t = q^{-1}$ gives $\Lambda(q\alk) \leq q\Lambda(\alk)$ and therefore
\[
 \Theta(\alpha_{k-1}) \leq q^{3/2} \Theta(\alk).
\]
Putting the last two equations together results into the bound
\[
C_\eta \| \Tdag\ek\|^2 + \alk\Delk \leq(2c_{\eta} C_1 q^3 + 8)\Theta(\alk)^2 =
(2c_{\eta} C_1 q^3 + 8)\alk\Lambda(\alk)^2.
\]
Firstly, this implies by omitting the second term on the left hand side that
\[
\| \Tdag\ek\|^2 \leq \dfrac{2c_{\eta} C_1 q^3 + 8}{C_\eta}\Theta(\alk)^2 \quad \text{and hence} \quad
C_1 \geq \dfrac{2c_{\eta} C_1 q^3 + 8}{C_\eta}.
\]
Hence, it is necessary that $C_\eta > 2q^3c_\eta$, which is equivalent to the inequality \eqref{eq:small_eta} assumed in the Lemma.
Then \eqref{eq:estimt} is true with $C_1 \geq \dfrac{8}{C_\eta-2q^3c_\eta}$.

Secondly, omitting the first term of the left hand side shows \hbox{$\Delk \leq(2c_{\eta} C_1 q^3 + 8)\Lambda(\alk)^2$}, so we have \eqref{eq:estimDelta} with $C_2 \geq 2c_{\eta} C_1 q^3 + 8$.
\bigskip\\
\emph{Case 2:} $\beta \alk\Delk^{1/2} \Lambda\paren{\frac{\| \Tdag\ek\|^2}{\Delk} + \alk}
\geq c_{\eta} \| \Tdag e_{\its-1}\|^2 + 4\Theta(\alk)^2$.\\
In this case \eqref{eq:induction:inequality} simplifies to 
\[
C_\eta \| \Tdag\ek\|^2 + \alk\Delk \leq 2 \beta \alk\Delk^{1/2}
\paren{\Lambda\left(\frac{\| \Tdag\ek\|^2}{\Delk} + \alk\right)}.
\]
Using again $\Lambda(0) = 0$ and the concavity we get $\Lambda(x) \geq \frac{x}{(a+b)}\Lambda(a+b)$ for all $0 \leq x \leq a + b$. Taking now $x = a$ and $x= b$ respectively implies $\Lambda(a) + \Lambda(b) \geq \Lambda(a+b)$. Thus we have
\begin{equation}\label{eq:secondcase}
C_\eta \| \Tdag\ek\|^2 + \alk\Delk \leq 2 \beta \alk\Delk^{1/2}
\paren{\Lambda\left(\frac{\| \Tdag\ek\|^2}{\Delk} \right) + \Lambda(\alk)}.
\end{equation}
It is again convenient to study two cases:\medskip\\
\emph{Case 2.1:} $\| \Tdag\ek\|^2 \leq \alk\Delk$.\\
Now the monotonicity of $\Lambda$ entails
\[
C_\eta \| \Tdag\ek\|^2 + \alk\Delk \leq 4 \beta \alk\Delk^{1/2} \Lambda\paren{\alk}.
\]
This shows that $\Delk^{1/2} \leq 4 \beta \Lambda\paren{\alk}$ and thereby \eqref{eq:estimDelta}
with $C_2 \geq 16\beta^2$. Plugging this into the right hand side of the last inequality and using the
case assumption for the left hand side we get
\[
(1+C_\eta)\| \Tdag\ek\|^2 \leq 16\beta^2 \alk \Lambda\paren{\alk}^2 = 16\beta^2 \Theta\paren{\alk}^2.
\]
Hence \eqref{eq:estimt} holds with $C_1 \geq \dfrac{16\beta^2}{1+C_\eta}$.\medskip\\
\emph{Case 2.2:} $\alk\Delk \leq \|\Tdag\ek\|^2$.\\
Dividing formula \eqref{eq:secondcase} by $\| \Tdag\ek\|$ results in
\[
C_\eta \| \Tdag\ek\| + \frac{\alk\Delk}{\| \Tdag\ek\|} \leq 2 \beta \alk
\left(\frac{\Delk}{\| \Tdag\ek\|^2}\right)^{1/2} \paren{\Lambda\left(\frac{\| \Tdag\ek\|^2}{\Delk}\right) + \Lambda(\alk)}.
\]
Since the functions $t^{-1/2}\Lambda(t)$ and $t^{-1/2}$ are monotonically decreasing, we obtain
\[
C_\eta \| \Tdag\ek\| + \frac{\alk\Delk}{\| \Tdag\ek\|} \leq 4 \beta \alk^{1/2} \Lambda(\alk).
\]
This shows that $C_\eta \| \Tdag\ek\| \leq 4 \beta \Theta(\alk)$, so \eqref{eq:estimt} is true with
$C_1 \geq \dfrac{16\beta^2}{C_\eta^2}$. Plugging this into the left hand side of the last equation gives
\[
\frac{\alk\Delk C_\eta}{4 \beta \alk^{1/2}\Lambda(\alk)} \leq 4 \beta \alk^{1/2} \Lambda(\alk).
\]
Now we see that $\Delk \leq 16\beta^2 \Lambda(\alk)^2/C_\eta$ and therefore that \eqref{eq:estimDelta} is valid with
$C_2 \geq 16\beta^2/C_\eta$. This completes the proof.\\

Now Theorem \ref{theo:main} follows easily:
\proof[Proof of Theorem \ref{theo:main}.]
The constant $C$ in the last lemma is independent of $\delta$ and $\gamma$. So if $\delta$ and $\gamma$ converge to $0$ in probability and if the probability that the tangential cone condition is not fulfilled goes to $0$, this implies convergence in probability of $\Delta(\widehat{\unknown}_{\capit},\xdag)$. That is the assertion of Theorem \ref{theo:main}.

\nocite{ChePou:12}
\bibliography{lit}
\bibliographystyle{humanbio}
\end{document}